\documentstyle[11pt,leqno]{article}
\newtheorem{th}{Theorem}[section]
\newtheorem{lem}[th]{Lemma}

\title{\Large\bf{Locally conformal flat Riemannian manifolds \\
with constant principal Ricci curvatures \\
and locally conformal  flat $\cal C$-spaces}}
\author{{\sc Stefan Ivanov} \thanks{The author supported
by Contract
MM 423/1994 with the Ministry of Science and Education of
Bulgaria and by
Contract 219/1994 with the University of Sofia "St. Kl.
Ohridski".} \hspace{5mm} {\sc Irina Petrova} \thanks{The author supported
by Contract
MM 413/1994 with the Ministry of Science and Education of
Bulgaria}}
\date{}
\begin{document}
\maketitle
\thispagestyle{empty}
\vspace{2mm}
\begin{center}
\end{center}
\vspace{25mm}
\begin{abstract}
It is proved that every locally conformal flat Riemannian manifold all 
of whose
 
Jacobi operators have constant eigenvalues along every geodesic
is with constant
principal Ricci curvatures.
A local classification (up to an isometry) of locally conformal flat
Riemannian manifold with constant Ricci eigenvalues is given in dimensions
$4,5,6,7$ and $8$.
It is shown that any $n$-dimensional $(4\leq n \leq 8)$ locally
conformal flat
Riemannian manifold with constant principal Ricci curvatures  is  a Riemannian
locally symmetric space.
\\[20mm]
{\bf Running title:}  Constant Ricci eigenvalues and ${\cal C}$-spaces.
\\[5mm]
{\bf Keywords:} Locally conformal
flat Riemannian manifolds, Constant eigenvalues of the Ricci tensor, Constant
eigenvalues of the Jacobi operator,
Curvature homogeneous spaces, Globally Osserman manifolds.
\\[5mm]
${\bf 1991 \quad MS \quad Classification: } 53B20; 53C15; 53C55; 53B35$
\end{abstract}
\newpage

\section{Introduction}

Curvature is a fundamental notion in Riemannian geometry. The Jacobi
operator is an important tool for studying the curvature.
If the Jacobi operator has  constant  eigenvalues then
the Riemannian manifold is
said to be globally  Osserman  manifold. The Osserman
conjecture states (see \cite{O}) that
any globally Osserman space is rank-one symmetric space. This conjecture has
been proved to be true for odd dimension and for dimensions two, four and
$2(2k+1), k \in N$ in \cite{Ch1}.
\par
J.Berndt and L.Vanhecke considered Riemannian manifolds
satisfying weaker  conditions in \cite{B-V1}. They introduced in \cite{B-V1}
the so-called $\cal C$-space
as a Riemannian manifold for which the Jacobi operators have constant
eigenvalues along every geodesic and the so-called $\cal P$-space
as a Riemannian manifold for which the Jacobi operators have parallel
eigenspaces along every geodesic.
The $\cal C$-spaces and $\cal P$-spaces can be regarded as a natural
generalization of Riemannian locally symmetric spaces since a Riemannian
manifold is locally symmetric iff it is a $\cal C$-space and it is a $\cal
P$-space simultaneously. There is a lot of examples of $\cal C$-spaces, namely
every naturally reductive Riemannian
homogeneous space as well as every commutative space is a $\cal
C$-space.  A local classifications (up to an isometry) of $\cal C$- spaces
and $\cal P$-spaces in dimensions $2$ and $3$ is given also in \cite{B-V1}.
\par
Globally Osserman spaces, $\cal C$-spaces, $\cal P$-spaces are of special
interest in the last
years. These spaces are studied in \cite{B-V1, B-V2, B-V3, B-P-V, Ch1, Ch2}.
\par
In this paper we consider locally conformal flat $\cal C$-spaces. Our main
observation is that these spaces have constant principal Ricci curvatures. We
give a local classification of
locally conformal flat Riemannian manifolds with constant principal Ricci
curvatures in  dimensions $4$,$5$,$6$,$7$ and $8$. Consequently, we obtain a
local description
of locally conformal flat $\cal C$-spaces in dimensions $4,5,6,7$ and $8$.
The aim of the paper is to prove the following
\begin{th}
\label{e11}
For an $n$-dimensional $(4\leq n\leq8)$ connected locally conformal flat
Riemannian manifold the following conditions are equivalent:
\par a) It is a ${\cal C}$-space;
\par b) It has constant principal Ricci curvatures;
\par c) It is locally (almost everywhere) isometric to one of the following
spaces:
\par\qquad i) a real space form;
\par\qquad ii) a Riemannian product of $1$-dimensional space and
of a real space form of dimen-
\par\qquad \quad sion $(n-1)$;
\par\qquad iii) a Riemannian product of two real space forms with
opposite constant sectional curvatures;
\par d) It is a Riemannian  locally symmetric space.
\end{th}
According to I.M.Singer \cite{Si} a Riemannian manifold $({\bf M},g)$ is called
{\it curvature homogeneous} if for any pair of points $p,q \in {\bf M}$ there
is a linear isometry $F : T_p{\bf M} \longrightarrow T_q{\bf M}$ between the
corresponding tangent spaces such that $F^* R_q = R_p$ (where $R$ denotes the
curvature tensor of type $(0,4)$).
We note that a locally conformal flat Riemannian manifold is curvature
homogeneous iff it has constant principal Ricci curvatures. Thus, every locally
conformal flat Riemannian manifold of dimension $4,5,6,7$ or $8$ is curvature
homogeneous iff it is a Riemannian locally symmetric space.
\par
Riemannian 3-manifolds  with constant
principal Ricci curvatures are studied in \cite{Ko,Ko1,KP,KP1,ST}. These
spaces are exactly the curvature homogeneous Riemannian 3-manifolds.
\par
Theorem \ref{e11} shows, unfortunately, that there are not interesting examples
of locally conformal flat $\cal C$-spaces in dimensions $4,5,6,7$ and $8$. But,
the situation with $\cal P$-spaces is completely different. In the recent
work \cite{ip3}, it is shown that there are exactly nine kinds of examples of
$4$-dimensional locally conformal flat $\cal P$-spaces and a lot
of these examples are not even curvature homogeneous.

\section {Locally conformal flat $\cal C$-spaces}

In this section we recall some definitions and prove our main
observation.
\par
Let $({\bf M},g)$ be an $n$-dimensional Riemannian manifold and $\nabla$
the Levi-Civita connection of the metric g. The curvature $R$ of $\nabla$
is defined by  $R(X,Y)=[\nabla_X,\nabla_Y]-\nabla_{[X,Y]}$ for
every vector fields $X,Y$ on {\bf M}. We denote
by $T_p{\bf M}$ the tangential space at a point $p\in{\bf M}$.
\par
Let $x\in
T_p{\bf M}$. The Jacobi operator is defined by
$$  \lambda_x(y)=R(y,x)x, \qquad y\in T_p{\bf M}.  $$
Let $\gamma(t)$ be a geodesic on ${\bf M}$ and $\dot{\gamma}$ denotes its
tangent vector field. We consider the family of smooth self-adjoint Jacobi
operators along $\gamma$ defined by $\lambda_{\dot{\gamma}}=R(X,
\dot{\gamma})\dot{\gamma}$ for every smooth vector field $X$ along $\gamma$.
A Riemannian manifold $({\bf M},g)$ is said to be a $\cal C$-space if the
operators $\lambda_{\dot{\gamma}}$ have constant eigenvalues along every
geodesic on ${\bf M}$.
\par
A Riemannian manifold $({\bf M},g)$ is said to be locally conformal flat
if around every point $p\in{\bf M}$ there exists a metric $\bar g$ which is
conformal to $g$ and $\bar g$ is flat. By the Weyl theorem, an $n$-
dimensional $(n\ge4)$ Riemannian manifold is locally conformal flat iff the
Weyl conformal tensor is zero (see e.g. \cite{Sch}) i.e.
the curvature tensor has the following form
\begin{equation}
\label{l1}
    R(x,y,z,u)=
-\frac{s}{(n-1)(n-2)}\left(g\left(y,z)g(x,u)-g(x,z)g(y,u\right)\right)+
\end{equation}
$$
\frac{1}{n-2}\left(Ric
\left(y,z)g(x,u)-Ric(x,z)g(y,u)+g(y,z)Ric(x,u)-g(x,z)Ric(y,u\right)\right),$$
$x,y,z,u\in T_p{\bf M},\quad p\in{\bf M}$,
where $Ric$ and $s$ are the Ricci tensor and the scalar curvature of $g$,
respectively.\\
The following condition also holds
\begin{equation}
\label{l2}
    (\nabla_xRic)(y,z)-(\nabla_yRic)(x,z)=
    \frac{1}{2(n-1)}(x(s)g(y,z)-y(s)g(x,z)).
\end{equation}
A Riemannian 3-manifold is locally conformal flat iff the condition (\ref{l2})
holds.
\par
The Ricci operator $\rho$ is defined by $g(\rho(x),y)=
Ric(x,y)$, $x$,$y\in T_p{\bf M}$, $p\in{\bf M}$. In every point
$p\in {\bf M}$ we consider
the Ricci operator as a linear self-adjoint operator on $T_p{\bf M}$. Let
$\Omega$ be the subset of ${\bf M}$ on which the number of distinct
eigenvalues of $\rho$ is locally constant. This set is open and dense in
${\bf M}$.  We can choose smooth eigenvalue functions
of $\rho$ on $\Omega$,
say $r_1,\ldots,r_n$, such that they form at each point of $\Omega$ the
spectrum of $\rho$  (see e.g. \cite{K, B-V1}). We fix $p\in\Omega$. Then
there exists a local orthonormal frame field $E_1,\ldots,E_n$ on an open
connected neighborhood $U$ of $p$ such that
\begin{equation}
\label{l3}
    \rho(E_i)=r_iE_i, \qquad i=1,2,\ldots,n.
\end{equation}
Our further considerations will take place in the neighborhood $U$.\\
For every $i,j,k,l \in \{1,2,...,n\}$, we set
$$  R_{ijkl}=R(E_i,E_j,E_k,E_l),\quad
    \omega_{ij}^{k}=g(\nabla_{E_i}E_j,E_k)  $$
$$  {(\nabla_iRic)}_{jk}=(\nabla_{E_i}Ric)(E_j,E_k).  $$
We have
\begin{th}\label{e21}
Every connected $n$-dimensional $(n\ge3)$ locally conformal flat
${\cal C}$-space has constant Ricci eigenvalues.
\end{th}
{\it Proof.} Let $p\in U$, $x\in T_p{\bf M}$ and $\gamma$ be a geodesic in
$U$ determined by the conditions $\gamma(0)=p$, ${\dot \gamma}(0)=x$.
Using (\ref{l1}), we get $trace(\lambda_{\dot \gamma})=Ric({\dot \gamma},
{\dot \gamma})$.       We have
$(\nabla_{\dot \gamma}Ric)({\dot \gamma},{\dot \gamma})=0$ along $\gamma$,
since $\bf M$ is a $\cal C$-space. At the point $p$, we obtain
\begin{equation}
\label{l4}
    (\nabla_xRic)(x,x)=0.
\end{equation}
The equality (\ref{l4}) holds
in every point $p\in U$ and for every tangent vector $x\in T_p{\bf M}$. We
get from (\ref{l4}) by a polarization that
\begin{equation}
\label{l5}
    (\nabla_iRic)_{jj}+2(\nabla_jRic)_{ji}=0,
    \qquad i,j\in\{1,\ldots,n\},\quad i\not=j.
\end{equation}
We obtain from (\ref{l2}) that
\begin{equation}
\label{l6}
    (\nabla_lRic)_{kk}+(\nabla_lRic)_{jj}-\frac{1}{n-1}E_l(s)-
    (\nabla_kRic)_{kl}-(\nabla_jRic)_{jl}=0,
\end{equation}
$$  j,k,l\in\{1,\ldots,n\},\quad j\not=k\not=l\not=j.  $$
Using (\ref{l5}) and (\ref{l6}), we derive the following equalities:
\begin{equation}
\label{l7}
    3(n-1)((\nabla_lRic)_{kk}+(\nabla_lRic)_{jj})-2E_l(s)=0,
\end{equation}
$$  j,k,l\in\{1,\ldots,n\},\quad j\not=k\not=l\not=j.  $$
We get from (\ref{l7}) that $(\nabla_lRic)_{kk}=(\nabla_lRic)_{jj}$,\quad
$j,k,l\in\{1,\ldots,n\}$, $j\not=k\not=l\not=j$. The latter equality implies
$E_l(r_k)=E_l(r_j)$,\quad $j,k,l\in\{1,\ldots,n\}$, $j\not=k\not=l\not=j$.
We obtain from (\ref{l4}) that $E_l(r_l)=0$,\quad $l\in\{1,\ldots,n\}$.
Substituting the last two equations into (\ref{l7}) we get $E_l(r_k)=0$,
\quad $l,k\in\{1,\ldots,n\}$ which completes the proof of the theorem.
\hfill {\bf Q.E.D.}

\section{Proof of Theorem 1.1}

We begin with the following technical results
\begin{lem}\label{l3.1}
Let $({\bf M},g)$ be an $n$-dimensional $(n\ge4)$ locally conformal flat
Riemannian manifold with constant Ricci eigenvalues such that at
least two of
them are distinct. Let $p\in{\bf M}$ and $U$ be the neighborhood
of $p$ described in the previous section.\\
Let $I_m$, $m=1,2,\ldots,l$, $(l\le n)$ be  subsets of the set
$\{1,\ldots,n\}$ such that
$$ \indent i)\qquad \bigcup_{m=1}^{l}I_m = \{1,\ldots,n\};
\hspace{9cm}$$ \par ii) $r_i=r_j$,\quad $i,j\in I_m$;
\par iii) $r_i\not=r_j$,\quad $i\in I_{m_1}$, $j\in I_{m_2}$, $m_1\not=m_2$.\\
Then, for any distinct $r,s,t\in\{1,\ldots,l\}$ we have
\begin{equation}
\label{l8}
    \omega_{ij}^{k}=0,\quad i,j\in I_s, k\in I_t;
\end{equation}
\begin{equation}
\label{l9}
    \omega_{ii}^{k}=0,\quad i\in I_s, k\in I_t;
\end{equation}
\begin{equation}
\label{l10}
    \omega_{im}^{k}=\frac{r_i-r_k}{r_m-r_k}\omega_{mi}^{k},
    \quad i\in I_s, m\in I_r, k\in I_t.
\end{equation}
\end{lem}
{\it Proof.} Using (\ref{l3}), we obtain from (\ref{l2}) that
\begin{equation}
\label{l11}
    (r_k-r_l)\omega_{ik}^{l}-(r_i-r_l)\omega_{ki}^{l}=0,
\end{equation}
\begin{equation}
\label{l12}
    (r_k-r_l)\omega_{kk}^{l}=0.
\end{equation}
The formulae (\ref{l8}), (\ref{l9}) and (\ref{l10}) follow from (\ref{l11}) and
(\ref{l12}) which proves the Lemma.\hfill {\bf Q.E.D.}\\
Further we have
\begin{lem}
\label{l3.2}
Let $({\bf M},g)$ be an $n$-dimensional $(n\ge4)$ locally conformal flat
Riemannian manifold with two distinct constant Ricci
eigenvalues. Then
$({\bf M},g)$ is locally isometric to one of the following spaces:
\par i)  ${\bf M^{n-1}\times M^1}$, where ${\bf M^{n-1}}$ is an
$(n-1)$-
dimensional real space form and ${\bf M^1}$ is an $1$-dimensional space;
\par ii) ${\bf M^m\times M^{n-m}}$, $(m<n)$ where ${\bf M^m}$ is an $m$-
dimensional real space form of constant sectional curvature ${\cal K}$ and
${\bf M^{n-m}}$ is an $(n-m)$-dimensional real space form of constant sectional
curvature ${\cal (-K)}$.
\end{lem}
{\it Proof.} Let $p\in{\bf M}$ and $U$ be as in Lemma \ref{l3.1}. \\
Let $r_1=r_2=\ldots=r_m\not=r_{m+1}=\ldots=r_n$. We set
$$  I_1=\{1,\ldots,m\},\quad I_2=\{m+1,\ldots,n\},  $$
$$  {\cal F}_1=span{\{E_i\}}_{i\in I_1},\quad
    {\cal F}_2=span{\{E_j\}}_{j\in I_2}.  $$
Each one of the smooth distributions ${\cal F}_1$ and ${\cal
F}_2$ is
autoparallel by Lemma \ref{l3.1}. Let $U^m$ and $U^{n-m}$ be the
correspondent integral submanifolds. The assertion of the Lemma follows
since $\bf M$ is locally conformal flat \hfill {\bf Q.E.D.}
\begin{lem}
\label{l3.3}
Let $({\bf M},g)$ be an $n$-dimensional $(n\ge4)$ locally conformal flat
Riemannian manifold with constant Ricci eigenvalues such that
exactly $l$ of them are distinct. Let $p\in{\bf M}$
and $U$ be a neighborhood of $p$ as in Lemma \ref{l3.1}.
Then, we have on $U$
\par i) $l\not=3$;
\par ii) there exists exactly one Ricci eigenvalue $r_s$ of
multiplicity $m_s$ such that $m_s>\frac{n}{2}$;
\par iii) $l=1,2$\quad or\quad $4\le l\le{n\over2}$\quad if $n$ is even;
\par\qquad $l=1,2$\quad  or\quad $4\le l\le{{n+1}\over{2}}$\quad if $n$ is odd;
\par iv)  if $l\ge4$ then $n\ge7$.
\end{lem}
{\it Proof.} Let $l\ge3$ and the first $l$-numbers of $r_1,\ldots,r_n$
be distinct and of multiplicity $m_1,\ldots,m_l$, respectively.
\\
Setting
\begin{equation}
\label{l13}
    u_k=2r_k-\frac{s}{n-1},\quad k=1,2,\ldots,l
\end{equation}
we claim that
\begin{equation}
\label{l14}
    u_k\not=0,
\end{equation}
\begin{equation}
\label{l15}
    u_k\not=u_j,\quad k,j\in\{1,\ldots,l\}, k\not=j;
\end{equation}
\begin{equation}
\label{l16}
    n-m_k=2u_k\sum_{{j=1}\atop{j\not=k}}^{l}\frac{m_j}{u_k-u_j},
    \quad k=1,2,\ldots,l.
\end{equation}
Formula (\ref{l15}) follows immediately from (\ref{l13}).
\\
To prove (\ref{l14}), we shall use the notations in Lemma \ref{l3.1}.
Let $s\not=t$, $s,t\in\{1,\ldots,l\}$,\quad $i\in I_t$,\quad
$j\in I_s$. We calculate using Lemma \ref{l3.1} that
\begin{equation}
\label{l17}
  R_{ijji}=2\sum_{{p=1}\atop{p\not\in{I_s}\cup{I_t}}}^{n}\frac{r_j-r_p}
  {r_i-r_p}(\omega_{ij}^{p})^2 .
\end{equation}
We obtain from (\ref{l17}) consequently:
$$
\sum_{{j=1}\atop{j\not\in I_t}}^{n}\frac{R_{ijji}}{r_i-r_j}=
\sum_{{s=1}\atop{s\not=t}}^{l}\frac{R_{ijji}}{r_i-r_j}=
\sum_{{s=1}\atop{s\not=t}}^{l}\sum_{j\in
I_s}\frac{2(r_j-r_p)}{(r_i-r_p)(r_i-r_j)} (\omega _{ij}^p)^2 =
$$
$$
2\sum_{{s,r=1}\atop{s\not=t,r\not=t,s\not=r}}^{l}\sum_{j\in
I_s,p\in I_r}\frac{r_j-r_p}{(r_i-r_p)(r_i-r_j)}(\omega _{ij}^p)^2,
$$
$$
\sum_{{j=1}\atop{j\not\in I_t}}^{n}\frac{R_{ijji}}{r_i-r_j}=
\sum_{{p=1}\atop{p\not\in I_t}}^{n}\frac{R_{ippi}}{r_i-r_p}=
\sum_{{r=1}\atop{r\not=t}}^{l}\sum_{p\in I_r}\frac{R_{ippi}}{r_i-r_p}=
$$
$$
\sum_{{r=1}\atop{r\not=t}}^{l}\sum_{p\in I_r}\sum_{{j=1}\atop{j\not\in I_t\cup
I_r}}\frac{2(r_p-r_j)}{(r_i-r_j)(r_i-r_p)}(\omega _{ip}^j)^2=
$$
$$
2\sum_{{s,r=1}\atop{s\not=t,r\not=t,s\not=r}}^{l}\sum_{p\in
I_r,j\in I_s}\frac{r_p-r_j}{(r_i-r_j)(r_i-r_p)}(\omega _{ij}^p)^2.
$$
We get comparing the latter  equalities that
\begin{equation}
\label{l18}
    \sum_{{j=1}\atop{j\not\in I_t}}^{n}\frac{R_{ijji}}{r_i-r_j}=0.
\end{equation}
On the other hand, we calculate using (\ref{l1}) that
\begin{equation}
\label{l19}
    \sum_{{j=1}\atop{j\not\in I_t}}^{n}\frac{R_{ijji}}{r_i-r_j}=
    \frac{1}{n-1}\sum_{{j=1}\atop{j\not=i}}^{l}\left(\frac{m_j(r_i+r_j)}
    {r_i-r_j}-\frac{m_js}{(n-2)(r_i-r_j)}\right) .
\end{equation}
We get (\ref{l16}) from (\ref{l18}) and (\ref{l19}). Since
$l\ge3$, then (\ref{l14}) follows from (\ref{l16}). We get from (\ref{l16})
that
\begin{equation}
\label{l20}
  \sum_{k=1}^{l}(n-m_k)m_ku_k=0; \qquad
  \sum_{k=1}^{l}\frac{(n-m_k)m_k}{u_k}=0.
\end{equation}
To prove i) let $l=3$. Then the first equality of (\ref{l20}) implies that any
two of the real numbers $u_1,u_2,u_3$ have the same sign and the third one
has the opposite sign. Let $u_2u_3>0$. If we set $x_i=\frac{u_i}{u_1}$,\quad
$i=2,3$ then we obtain from (\ref{l20}) that
$$  (n-m_2)m_2x_2+(n-m_3)m_3x_3=-(n-m_1)m_1,  $$
$$  (n-m_2)m_2x_3+(n-m_3)m_3x_2=-(n-m_1)m_1x_2x_3.  $$
We get multiplying the latter two equalities that
\begin{equation}
\label{l21}
    (n-m_2)m_2(n-m_3)m_3(x_2-x_3)^2+(4m_1m_2m_3(n-m_1)+4m_2^2m_3^2)x_2x_3=0.
\end{equation}
The left hand side of (\ref{l21}) is strictly positive since
$x_2x_3>0$. This contradiction proves i).\\
We get from (\ref{l16}) that
\begin{equation}
\label{l22}
    n-m_k=2\sum_{{j=1}\atop{j\not=k}}^{l}\frac{m_ju_j}{u_j-u_k} .
\end{equation}
Using (\ref{l22}), we calculate
\begin{equation}
\label{l23}
    \sum_{k=1}^{l}(n-2m_k)m_ku_k^2=-{\left(\sum_{k=1}^{l}m_ku_k\right)}^2.
\end{equation}
Now, ii) follows immediately from (\ref{l23}). The assertion iii)
and iv) are
consequences of i) and ii). Thus, the whole Lemma is proved. \hfill {\bf
Q.E.D.}
\par
We are ready to prove the equivalence between b) and c).
\par
If the Ricci tensor has exactly one eigenvalue then $\bf M$
is an Einstein space. Hence, $\bf M$ is a
real space form since it is locally conformal flat.
\par
If the Ricci tensor has exactly two distinct eigenvalues then
the assertion follows from Lemma \ref{l3.2}.
\par
If $n=4,5,6$, then Lemma \ref{l3.3} implies that the Ricci tensor has
either one or two distinct eigenvalues and the assertion follows.
\par
Let $n=7$ or $8$. It follows from Lemma \ref{l3.3} that the Ricci
tensor has one eigenvalue, or two distinct eigenvalues, or four distinct
eigenvalues such that three of them are of the
multiplicity one. We shall prove that the latter case is
impossible for any $n\ge4$.\\
Let we assume $r_1\not=r_2\not=r_3\not=r_4=\ldots=r_n$. Then $m_4=n-3$. We
obtain from (\ref{l16}) that
$$  3=2u_4\sum_{i=1}^3\frac{1}{u_4-u_i}.  $$
The latter equality is equivalent to the equality
\begin{equation}
\label{l24}
    3u_4^3-u_4^2\sum_{i=1}^3u_i-u_4\sum_{{i,j=1}\atop{i<j}}^3u_iu_j+
    3u_1u_2u_3=0 .
\end{equation}
We get from (\ref{l20}) that
\begin{equation}
\label{l25}
    \sum_{i=1}^3u_i=-\frac{3(n-3)}{n-1}u_4,
\end{equation}
\begin{equation}
\label{l26}
    \sum_{i=1}^3\frac{1}{u_i}=-\frac{3(n-3)}{n-1}\frac{1}{u_4} .
\end{equation}
Using the latter two equalities, we obtain from (\ref{l24}) that
\begin{equation}
\label{l27}
    u_4^3=-u_1u_2u_3 .
\end{equation}
Setting $x_i=\frac{u_i}{u_4}$,\quad $i=1,2,3$, we derive from (\ref{l25}),
(\ref{l26}) and (\ref{l27}) that
\begin{equation}
\label{l28}
    \sum_{i=1}^3x_i=-\frac{3(n-3)}{n-1};\quad
    \sum_{{i,j=1}\atop{i<j}}^3x_ix_j=\frac{3(n-3)}{n-1};\quad
    x_1x_2x_3=-1.
\end{equation}
The equalities (\ref{l28}) imply that $x_1$, $x_2$, $x_3$ are
the roots of the following cubic equation
$$  x^3+\frac{3(n-3)}{n-1}x^2+\frac{3(n-3)}{n-1}x+1=0.  $$
It is easy to see that this cubic equation has exactly one real and two
complex roots. This contradiction proves the equivalence of b) and c).
\par
Now, it follows from Theorem \ref{e21} that
every $n$-dimensional $(4\le n\le8)$ locally conformal flat ${\cal C}$-space
has to be locally isometric to one of the spaces described in i), ii) and iii).
It is clear that every Riemannian manifold of the type i), ii) or iii) is
locally symmetric. Hence, it is a ${\cal C}$-space.
\par
Conversely, every locally conformal flat Riemannian locally symmetric
space is locally isometric to one of the spaces described by i), ii) or iii).
This completes the proof of Theorem 1.1.
\hfill \qquad {\bf Q.E.D.}
\par
{\bf Remark}. Riemannian manifolds for which the natural skew-symmetric
curvature operators $\kappa_{x,y}(z):=R(x,y)z, \quad x,y,z\in T_p{\bf M}, p\in
{\bf M}$ have constant eigenvalues along every unit circle are called $\cal
O$-spaces
(see \cite{ip1}).  Riemannian manifolds for which the natural skew-symmetric
curvature operators $\kappa_{x,y}$ have parallel Jordanian bases along every
unit circle are called $\cal T$-spaces (see \cite{ip2}). Following the proof of
Theorem \ref{e21}, it is not difficult to see that any $4$-dimensional locally
conformal flat $\cal O$-space as well as any $4$-dimensional locally conformal
flat $\cal T$-space has constant Ricci eigenvalues. Hence, any $4$-dimensional
locally
conformal flat $\cal O$-space as well as any $4$-dimensional locally conformal
flat $\cal T$ space is locally (almost everywhere) isometric to one of the
spaces described in c) of Theorem \ref{e11}.

\begin{flushleft}
{\bf Authors' address:}\\
         Stefan Ivanov,Irina Petrova,\\
         University of Sofia,Faculty of Math. and Inf.,\\
         bul. James Boucher 5, 1126 Sofia,\\
         BULGARIA\\
         E-mail:\quad ivanovsp@fmi.uni-sofia.bg\\
         E-mail:\quad ihp@vmei.acad.bg
\end{flushleft}


\begin{thebibliography}{9}

\bibitem{B-V1} J.Berndt and L.Vanhecke. {\em Two naturally  generalizations
of locally symmetric spaces.} Diff.Geom. and Appl.{\bf 2}(1992),
57-80.
\bibitem{B-V2} J.Berndt and L.Vanhecke. {\em Geodesic spheres and
generalizations of symmetric spaces.} Boll.Un. Mat.Ital.A(7), {\bf 7}
(1993), no.1, 125-134.
\bibitem{B-V3} J.Berndt and L.Vanhecke. {\em Geodesic sprays and $\cal C$-
and $\cal B$-spaces.} Rend.Sem. Politec.Torino {\bf 50}(1992), no.4,
343-358.
\bibitem{B-P-V} J.Berndt, F.Pr\"ufer and L.Vanhecke. {\em Symmetric-like
Riemannian manifolds and geodesic symmetries.} Preprint.
\bibitem{Ch1} Q.-S.Chi. {\em A curvature characterization of certain locally
rank-one symmetric spaces.} J.Diff.Geom.{\bf 28}(1988), 187-202.
\bibitem{Ch2} Q.-S.Chi. {\em Curvature characterization and  classification
of rank-one symmetric spaces.} Pacific J.Math., vol.{\bf 150}(1991), no.1, 
31-41.
\bibitem {ip1} S.Ivanov, I.Petrova, {\em Riemannian manifold in which certain
curvature operator has constant eigenvalues along each circle},
Ann. Glob. Ann. Geom., to appear.
\bibitem {ip2} S.Ivanov, I.Petrova, {\em Curvature operator with
parallel Jordanian basis on circles},  Riv. Mat. Univ. Parma, to appear.
\bibitem {ip3} S.Ivanov, I.Petrova, {\em Locally conformal flat
Einstein-like
4-manifolds  and locally conformal flat Riemannian 4-manifolds
all of whose  Jacobi
operators have parallel eigenspaces along every geodesic},
to appear.
\bibitem{K} T.Kato. Perturbation theory for linear operators. {\em
 Springer,Berlin,1966}.
\bibitem{Ko} O.Kowalski. {\em A classification of Riemannian $3$-manifolds
with constant principle Ricci curvatures\quad $\rho_1=\rho_2\not=\rho_3$}.
Nagoya Math.J.{\bf 132}(1993), 1-39.
\bibitem{Ko1} O.Kowalski {\em Nonhomogeneous Riemannian $3$-manifolds with
distinct constant Ricci eigenvalues}. Comment.Math.Univ.Carolinae
{\bf 34}(1993), 451-457.
\bibitem{KP} O.Kowalski and F.Pr\"ufer. {\em On Riemannian $3$-manifolds
with distinct Ricci eigenvalues}. Math.Ann.{\bf 300}(1994), 17-28.
\bibitem{KP1} O.Kowalski and F.Pr\"ufer. {\em A classification of special
Riemannian 3-manifolds with distinct constant Ricci eigenvalues}. J. for
Analysis and its Applications, vol.{\bf 13}(1994), no.4, 1-6.
\bibitem{O} R.Osserman. {\em Curvature in the eighties.} Amer.Math.Monthly
{\bf 97}(1990), 737-756.
\bibitem{Si} I.M.Singer, {\em Infinitesimally homogeneous
spaces}, Commun. Pure Appl. Math ${\bf 13} (1960), 685-697.$
\bibitem{Sch} J.A.Schouten.  {\em Ricci calculus} 2nd ed.,
Springer 1954.
\bibitem{ST} A.Spiro and F.Tricceri. {\em 3-dimensional Riemannian metric
with prescribed Ricci principal curvatures}. J.Math.Pures Appl.,
vol.{\bf 74}(1995), 253-271.
\end{thebibliography}
\end{document}